\documentstyle[11pt]{article}

\newcommand{\be}{\begin{equation}}
\newcommand{\ee}{\end{equation}}
\newcommand{\ba}{\begin{eqnarray}}
\newcommand{\ea}{\end{eqnarray}}
\newcommand{\baa}{\begin{eqnarray*}}
\newcommand{\eaa}{\end{eqnarray*}}
\newcommand{\bb}{\begin{thebibliography}}
\newcommand{\eb}{\end{thebibliography}}

\newcommand{\bi}[1]{\bibitem{#1}}
\newcommand{\lab}[1]{\label{#1}}
\newcommand{\re}[1]{(\ref{#1})}

\newcommand\bt{\beta}
\newcommand\om{\omega}
\newcommand\la{\lambda}

\begin{document}

\begin{flushright}
CRM-2278\\
March 1995
\end{flushright}

\begin{center}

\bigskip\bigskip

{\large \bf $q$-Ultraspherical Polynomials for $q$ a Root of Unity}%
\footnote{Published in {\bf Lett. Math. Phys. 37 (1996), pp. 173-180.}}

\bigskip\medskip

{\large \bf Vyacheslav Spiridonov%
\footnote{On leave from the Institute for Nuclear Research, Russian Academy
          of Sciences, Moscow, Russia}
}
\medskip

{\em Centre de Recherches Math\'ematiques, Universit\' e de Montr\' eal, \par
C.P. 6128, Succ. Centre-ville, Montr\' eal, Qu\' ebec, H3C 3J7, Canada \par
e-mail: spiridonov@lps.umontreal.ca}

\bigskip\medskip

{\large \bf Alexei Zhedanov}

\medskip                                                                        

{\em Physics Department, Donetsk University, Donetsk, 340055, Ukraine}

\end{center}

\begin{abstract}
Properties of the $q$-ultraspherical polynomials for $q$ being  a 
primitive root of unity are derived using a formalism of the $so_q(3)$ 
algebra. The orthogonality condition for these polynomials provides a 
new class of trigonometric identities representing discrete 
finite-dimensional analogs of $q$-beta integrals of Ramanujan.

\end{abstract}




{1. \it Introduction.}
There are many applications of quantum algebras, or $q$-analogues
of Lie algebras in mathematical physics \cite{kulish}. In particular,
it is well known that the corresponding representation theory 
is related to $q$-special functions and $q$-orthogonal polynomials. 
In this Letter we investigate from this point of view a class of
$q$-ultraspherical polynomials for $q$ a root of unity.

The monic $q$-ultraspherical polynomials satisfy three-term recurrence relation
\be
P_{n+1}(x)+u_nP_{n-1}(x)=xP_n(x), \quad n=1, 2, \dots
\lab{rec} \ee
and initial conditions $P_0(x)=1,\; P_1(x)=x,$
where coefficients $u_n$ have the form \cite {AI,AW, GR}
\be
u_n={{(1-q^n)(1-q^{n+2\beta-1})}\over{(1-q^{n+\beta})(1-q^{n+\bt-1})}}.
\lab{un} \ee

These polynomials are well investigated for $0\le q\le 1$. 
The case of $q$ a root of unity was considered only in the nontrivial
limits $\bt \to 0$ or 1 leading to the so-called sieved ultraspherical 
polynomials \cite{AAA}. In these limits the 
recurrence coefficients $u_n>0$ for all $n>0$, so that 
the sieved polynomials are orthogonal
on an interval of real line with some continuous measure. 
Consider what happens if 
$\bt \ne 0, 1$ and $q$ is $N$th primitive
root of unity, $q=\exp 2\pi ip/N, \; (p, N)=1$. 
Then $u_n$ can be rewritten in the form
\be
u_n={{\sin\om n \sin\om(n+2\bt-1)}\over{\sin\om(n+\bt)\sin\om(n+\bt-1)}}, 
\lab{un2} \ee
where $\om=\pi p/N.$ 
In order for $P_n(x)$ to possess a positive weight function
(i.e. to be ``classical" orthogonal polynomials) $u_n$ must be positive. 
Reality of $u_n$ implies reality of the parameter $\bt$. Then \re{un2} 
cannot satisfy the positivity condition for
all $n>0$. One has to analyze thus the finite-dimensional cases when
$u_0=u_M=0$ for some natural $M\geq 2$ and $u_n>0$ for $0<n<M$ (then $M$ is the
dimension of Jacobi matrix defining the eigenvalue problem \re{rec}). 

The possibility to construct finite-dimensional classical orthogonal 
polynomials on trigonometric
grids (including those generated by $q$ a root of unity) was mentioned in
\cite {NSU}, however no explicit examples were presented. 
For examples of non-standard 
$q$-special functions at $q^N=1$ related to the Lam\'e-type 
differential equations, see \cite{skorik}. 

Recently we have found that a special class of $q$-ultraspherical polynomials 
for $q^N=1$ can be obtained by ``undressing" of the finite-dimensional
Chebyshev polynomials \cite{ann,dar}, existence of the more complicated 
cases has been pointed out as well. Here we describe $q$-ultraspherical
polynomials with positive measure when $q$ is the $p=1$ primitive root of 
unity, i.e. when
\be
q=\exp 2\pi i/N,\qquad N=2, 3, \dots .
\lab{q} \ee
The main tool in our analysis will be the theory of 
representations of quantum algebra
$so_q(3)$ introduced in \cite{ode,fairlie}. The connection of this algebra
with some (finite-dimensional) class of $q$-ultraspherical polynomials
for real $q$ was mentioned in \cite{z}. 

{2. \it Restrictions upon the Parameters.}
Consider restrictions upon the parameters $M$
and $\bt$ necessary for positivity of $u_n$ in the interval $0<n<M$.
There are two principally different situations. 
The first possibility is $M=N,$ 
when dimension of representation coincides with $N$. Then the
positivity condition leads to the following restrictions upon $\bt$
(all the inequalities are modulo $N$):
\be
-1/2<\bt<0, \qquad 
0<\bt<1, \qquad
1<\bt<3/2.
\lab{bt+} \ee
Note that the ends of intervals do not belong to this class.
Indeed, there are ambiguities 
in choosing the coefficients $u_n$ for $\bt=0,1$. 
The cases $\bt=-1/2,\; 3/2$ lead to $u_M=0$ for $M=2, \; N-2,$
which is not allowed since we demand the dimension of 
representation to be equal to $N$. 

The second possibility arises if $2\bt$ is a natural number:
\be
2\bt=j= 1, 2, 3, \dots , N-1.
\lab{beta2} \ee
We define coefficients $u_n$ in the case $\bt=1$ as 
$u_0=u_{N-1}=0, \; u_1=\dots=u_{N-2}=1.$ 
Then the dimension of representations defined by (\ref{beta2})
is determined by the general formula $M=N+1-j.$ 
We note by passing that the even $j$ cases describe undressing of the
``discrete finite well" in discrete quantum mechanics \cite{zakh}. 
The case $j=1$ corresponds to the $q$-Legendre polynomials because 
$u_n$ reduce to the ordinary Legendre polynomials' recurrence coefficients
for $N\to\infty$. This system has an interesting physical application
in the Azbel-Hofstadter problem of electron on a two-dimensional lattice 
in a magnetic field \cite{wig}. 

It is not difficult to see that \re{bt+} and \re{beta2} are the only 
choices of $\beta$ providing positive $u_n$ (we skip the trivial case 
$\beta=-1/2$ when $M=2$ and assume $\beta\neq 0$). 

{3. \it  $so_q(3)$ Algebra and Its Representations.} 
Consider the associative 
algebra of three generators defined by commutation relations
\be
[K_0,K_1]_{\om} = K_2 ,\qquad
[K_1,K_2]_{\om}=-K_0 ,\qquad
[K_2,K_0]_{\om}=-K_1,
\lab{com} \ee
where $[A,B]_{\om}=e^{i\om /2}AB - e^{-i\om/2}BA$
denotes so-called $q$-mutator and $\om=\pi/N$.
Note that if $K_0$ and $K_1$ are
hermitian operators then $K_2$ should be anti-hermitian. Note also that for
$N\to\infty$ we get the ordinary rotation algebra $so(3)$, which justifies the
name $so_q(3)$ for the algebra \re{com}. This algebra was introduced
(in a slightly different form) in 1986 by Odesskii \cite{ode} who considered
its representations without discussion of unitarity. 
Its relations to quantum algebra
$su_q(2)$ and $q$-special functions have been discussed in \cite{fairlie,z}.
A non-compact version of this algebra $so_q(2, 1)$ for real $q$ appeared
to be quite useful as an algebra describing dynamical symmetry of
special classes of discrete reflectionless potentials \cite{ann}. 

Construct the unitary finite-dimensional representations of $so_q(3)$ by
the method, proposed in \cite{z}.
It is easy to see that there exists an orthonormal basis $|n\rangle$ for 
which the operator 
$K_0$ is diagonal whereas the operator $K_1$ is two-diagonal:
\be
K_0|n\rangle = \la_n|n\rangle, \qquad K_1|n\rangle =a_{n+1}|n+1\rangle + 
a_n|n-1\rangle,
\lab{rep} \ee
where $\la_n$ are eigenvalues of operator $K_0$ and $a_n$ are matrix
elements of a representation.
Substituting (\ref{rep}) into (\ref{com}) 
we find 
\be
\la_n={{\cos\om(n+\bt)}\over{\sin\om}}, \quad n=0, 1, \dots, M-1,
\lab{lan} \ee
\be
a_n^2={{\nu - \cos\om(n+\bt)\cos\om(n+\bt-1))}\over{4\sin^2\om
\sin\om(n+\bt)\sin\om(n+\bt-1)}},
\lab{an} \ee
where $\nu$ is the eigenvalue of Casimir operator of $so_q(3)$ algebra
which can be represented in the form
\be 
Q=\frac{1}{2}(K_2\tilde K_2 + \tilde K_2 K_2)-\cos \omega (K_0^2+K_1^2),
\lab{cas}\ee
where $\tilde K_2=[K_0, K_1]_{-\omega}$.
The condition $a_0=0$ gives 
$\nu= \cos\om\bt\cos\om(\bt-1)$, and then
\be
a_n^2={{\sin\om n\sin\om (n+2\bt-1)}\over{4\sin^2\om \sin\om (n+\bt)
\sin\om (n+\bt-1)}}
\lab{an2} \ee
coincide up to a constant factor with \re{un2}.
Hence, all unitary irreducible finite-dimensional 
representations of the $so_q(3)$ algebra for $q=\exp{2\pi i/N}$ 
are exhausted by two possibilities: \re{bt+} of the 
dimension $N$ and \re{beta2} of a smaller dimension.

We introduce for brevity the following terminology: the representations
\re{bt+} will be called ``complementary series"; 
the cases \re{beta2} will be called ``integer series"
for $\bt=1,2,3,\dots $ and ``half-integer series" for $\bt=1/2,3/2,\dots $
(the case $\bt=1/2$ is exceptional, for it belongs to both
complementary and half-integer series). 
Such definitions are supported by the fact that eigenvalues of the 
Casimir operator (\ref{cas}) are continuous and quantized for the 
complementary and integer (half-integer) series respectively
(cf. with the Lorentz algebra).

{4. \it $q$-Ultraspherical Polynomials as Overlap Coefficients.}
$q$-Ultraspherical polynomials arise as 
overlap coefficients between two distinct bases in the space of 
$so_q(3)$ representations.
Indeed, because $K_1$ is hermitian, we can choose another orthonormal basis
$|s)$ for which the operator $K_1$ is diagonal:
\be
K_1 |s)=\mu_s |s).
\lab{mu} \ee
$so_q(3)$ is a specisal case of the $AW(3)$ algebra describing 
symmetries of the Askey-Wilson polynomials \cite{zaw}. 
From general properties of the latter algebra we know that 
the operator $K_0$ cannot be more than tridiagonal in this dual basis:
\be
K_0|s)= d_{s+1}|s+1) + d_s|s-1) + b_s|s),
\lab {K0} \ee
where $s=0,\dots ,M-1$. 
The explicit form of matrix elements $d_s$ and $b_s$ depends on the
representation series of $so_q(3)$. 
For integer and half-integer series we have the spectrum
\be
\mu_s= {{\cos\om (s+j/2)}\over{\sin\om}},
\lab{mus} \ee
where $j=2\bt=1, 2, 3, \dots,$ and the matrix elements
\be
d_s^2={{\sin\om s \sin\om(s+j-1)}\over{4\sin^2\om \sin\om(s+j/2)
\sin\om(s-1+j/2)}} , \qquad b_s\equiv 0.
\lab {selem} \ee
This case is symmetric -- matrix elements in the bases $|n\rangle$
and $|s)$ are identical. 

For the complementary series we have
\be
\mu_s={{\cos\om(s+1/2)}\over{\sin\om}},\qquad s=0,1,\dots, N-1, 
\lab{mus2} \ee
\be
d_s^2={{\sin\om(s-\bt+1/2)\sin\om(s+\bt-1/2)}\over{4\sin^2\om \sin\om(s+1/2)
\sin\om(s-1/2)}},\qquad s=1, 2, \dots, N-2, 
\ee
\be
b_0=b_{N-1}={{\sin\om(1/2-\bt)}\over{2\sin\om \sin\om/2}}, \qquad 
d_0=d_{N-1}=b_1=b_2=\dots=b_{N-2}=0.
\lab{selem2} \ee
Note that the coefficients $d_n$ do not truncate automatically at the
ends of the index intervals, one has to do it by hands.
In this case the operator $K_0$ is {\it tridiagonal} instead of being 
two-diagonal as it is expected from the symmetry of $so_q(3)$ 
(after the replacement $K_2\to iK_2$ this algebra looks totally symmetric
under the cyclic permutations). This anomaly shows that permutations of 
generators are not necessarily unitary automorphisms
of the ``Cartesian" quantum algebras.

For any series: integer, half-integer or the complementary one, we can find the
overlap coefficients between two bases $(s|n\rangle$. These coefficients can 
be factorized, $(s|n\rangle =(s|0\rangle S_n(\mu_s),$ 
where $S_n(\mu_s)$ are symmetric polynomials satisfying the
three-term recurrence relation
\be
a_{n+1}S_{n+1} + a_nS_{n-1} = \mu_sS_n, \qquad S_0=1,\quad S_1=\mu_s/a_1.
\lab{recs} \ee
Polynomials $S_n(\mu_s), \; n>1$ are connected with $P_n(x)$ \re{rec}
by the simple relation 
\be
P_n(x)=\sqrt{u_1u_2\cdots u_n}S_n(\mu_s), \qquad x=2\mu_s \sin\omega.
\ee

Existence of the algebraic interpretation of $q$-ultraspherical 
polynomials in terms of the $so_q(3)$ algebra allows to
calculate the weight function for these polynomials. Omitting the details
(we use the method described, e.g. in \cite{zaw}), we present 
the final result.

For the integer and half-integer series the weight function has the form
\be
w_s(j)=\sin\om(s+j/2)\prod_{l=1}^{j-1}\sin\om(s+l), \quad
s=0, 1, \dots,  N-j.
\lab{wint} \ee
For the complementary series we have
\be
w_s(2\beta) = 
w_0{{\sin\om(s+1/2)}\over{\sin\om/2}}\prod_{l=0}^{s-1} {{\sin\om(\bt
+1/2+l)}\over{\sin\om(-\bt+3/2+l)}},\quad s=1, 2, \dots  , N-1,
\lab{wcomp} \ee
where $w_0$ is (undefined) 
value of the weight function at the points $s=0$ and $s=N-1$.

Clearly both weight functions \re{wint} and \re{wcomp} are known up to
a normalization constant. This constant can not be found directly from the
representation theory of $so_q(3)$ algebra and should be calculated 
separately. For the integer and half-integer series
it is possible to derive weight functions together with
normalization constants with the help of Darboux transformations. 

{5. \it Darboux Transformation and Normalization Constants.}
Let $P_n^{(j)}(x)$ be the 
monic $q$-ultraspherical polynomials belonging to integer
or half-integer series. Write the orthogonality relation in the form
\be
\sum_{s=0}^{N-j}P_n^{(j)}(x_s)P_m^{(j)}(x_s)w_s(j) = h_n(j)\delta_{nm},
\lab{orth} \ee
where $h_n(j)$ are normalization constants which have to be found, and 
$w_s(j)$ is
given by \re{wint}. Recall that for the taken series
$x_s(j)=2\cos\om(s+j/2).$ 

The crucial observation is that the polynomials $P_n^{(j+2)}(x)$ 
and $P_n^{(j)}(x)$ are related to each other by the Darboux 
transformation \cite{dar} which is equivalent in our case to the 
transition to symmetric kernel polynomials \cite{chi}:
\be
P_n^{(j+2)}(x_s(j+2))={{P_{n+2}^{(j)}(x_s(j))-A_n(j)P_n^{(j)}(x_s(j))}\over
{x_s^2(j)-x_0^2(j)}},
\lab{darbu} \ee
where
\be
A_n(j)={{P_{n+2}^{(j)}(x_0(j))}\over{P_n^{(j)}(x_0(j))}}.
\ee
From the theory of kernel polynomials \cite{chi} it follows that
the weight function is transformed as
\be
w_s(j+2)=w_{s+1}(j)(x_0^2(j)-x_{s+1}^2(j))/4,
\lab{trw} \ee
and the normalization constants are transformed as
\be
h_n(j+2)=h_n(j)A_n(j)/4.
\lab{nconst} \ee

In our case the factor $A_n(j)$ can be easily found
\be
A_n(j)={{\sin\om(n+j+1) \sin\om(n+j)}\over
{\sin\om(n+1+j/2) \sin\om(n+j/2)}}.
\lab{An} \ee
Then using \re{An}, \re{nconst} and  starting from $j=2$ 
and $j=1$ we obtain explicit expressions
for the normalization constant in the case of integer series, $j=2k$: 
\be
h_n(2k)={{h_0(2)s_{n+k+1}s_{n+k+2}\dots s_{n+2k-1}}\over
{4^{k-1}s_{n+1}s_{n+2}\dots s_{n+k-1}}}, \qquad n=0, 1, \dots , N-2k,
\lab{normint} \ee
where $s_n\equiv\sin\om n$,
and in the case of half-integer series, $j=2k+1$:
\be
h_n(2k+1)={{h_0(1)s_1^2s_2^2\dots s_n^2s_{n+1}s_{n+2}\dots s_{n+2k}}\over
{4^ks_{1/2}s_{3/2}^2s_{5/2}^2\dots s_{n+k-1/2}^2s_{n+k+1/2}}}, \quad
n=0, 1, \dots, N-2k-1,
\lab{normhalf} \ee
where we assume that for $n=0$ the product $s_1^2s_2^2\dots s_n^2$ is 
replaced by 1.
 The only coefficients left to be determined are $h_0(2)$ and $h_0(1)$. But
these constants can be calculated directly. Indeed, for $j=2$ we
have the finite-dimensional Chebyshev polynomials
\be
P_n(x_s)={{\sin\om(n+1)(s+1)}\over{\sin\om(s+1)}}.
\ee
The weight function in this case is $w_s(2)=\sin^2\om(s+1).$ Hence 
\be
h_0(2)=\sum_{s=0}^{N-2} \sin^2\om(s+1) = N/2.
\lab{h(2)} \ee

The $q$-Legendre polynomials case $j=1$ is characterized by the surprisingly 
simple weight function, $w_s(1)=\sin\om(s+1/2),$ from which we derive 
\be
h_0(1)=\sum_{s=0}^{N-1}\sin\om(s+1/2)={{1}\over{\sin\om/2}}.
\lab{h(1)} \ee
Thus we have calculated $h_n(j)$ for both integer 
\re{normint} and half-integer \re{normhalf} series.

The $n=0$ values of the normalization constants $h_n(j)$ provide
non-trivial trigonometric identities of the form 
$\sum_{r=0}^{M-1}w_r(j)=h_0(j).$ For $j=2k$ we get 
\begin{equation}
\sum_{r=0}^{N-2k}s_{r+k}s_{r+1 }s_{r+2}\dots s_{r+2k-1}=
{2Ns_{k+1}s_{k+2}\dots s_{2k-1}\over 4^k s_1s_2\dots s_{k-1}},
\label{id1}
\end{equation}
where $k=2, 3, \dots, [N/2]$. For $j=2k+1$ we get 
\begin{equation}
\sum_{r=0}^{N-2k-1}s_{r+k+1/2}s_{r+1 }s_{r+2}\dots s_{r+2k}=
{s_1s_2\dots s_{2k} \over 4^k s_{1/2}^2s_{3/2}^2\dots s_{k-1/2}^2s_{k+1/2}},
\label{id2}
\end{equation}
where $k=1, 2, \dots, [N/2].$
These identities can be considered as discrete finite-dimensional 
analogs of the $q$-beta integrals of Ramanujan 
(the latter are defined for real $q$, see e.g. \cite{AW}). 
Of course, in the case of $q$-ultraspherical
polynomials we get only special cases of such identities.
Their generalizations are obtained if we deal with general 
trigonometric (four-parameter) analogs of the Askey-Wilson polynomials. 
Existence of similar simple relations for the representations of complementary 
series is an interesting open problem which we hope to address in the future. 

{\it Acknowledgments.}
The authors are indebted to R.Askey, D.Masson, M.Rahman, S.Suslov, 
and L.Vinet for stimulating discussions.
The work of V.S. is supported by NSERC of Canada and by Fonds FCAR
of Qu\'ebec.
The work of A.Zh. is supported in part by the ISF Grant No. U9E000.


\bb{99}

\bi{kulish} {\it Quantum Groups}, Lect. Notes in Math. v. 1510, 
ed. P.P.Kulish (Springer-Verlag, Berlin, 1992).

\bi{AI} Askey R. and Ismail M., {\it A generalization of ultraspherical 
polynomials}, Studies in Pure Mathematics, ed. P.Erd\"os 
(Birkh\"auser, Basel, 1983) pp. 55-78. 

\bi{AW} Askey R. and Wilson J., {\it Some basic hypergeometric orthogonal
polynomials that generalize Jacobi polynomials}, Mem. Am. Math. Soc.
{\bf 54}, 1-55 (1985).

\bi{GR} Gasper G. and Rahman M., {\it Basic Hypergeometric Series}
(Cambridge University Press, 1990).

\bi{AAA} Al-Salam W., Allaway W., and Askey R., 
{\it Sieved ultraspherical polynomials}, 
Trans. Amer. Math. Soc. {\bf 284}, 39-55  (1984). 

\bi{NSU} Nikiforov A.F., Suslov S.K., and Uvarov V.B., {\it Classical Orthogonal
Polynomials of Discrete Variable} (Springer-Verlag, Berlin, 1991).

\bi{skorik} Skorik S. and Spiridonov V., {\it Self-similar potentials and
the $q$-oscillator algebra at roots of unity}, Lett. Math. Phys.
{\bf 28}, 59-74 (1993).

\bi{ann} Spiridonov V. and Zhedanov A., {\it Discrete reflectionless
potentials, quantum algebras, and $q$-orthogonal polynomials},
Ann. Phys. (N.Y.) {\bf 237}, 126-146 (1995).

\bi{dar} Spiridonov V. and Zhedanov A., {\it Discrete Darboux transformation,
the discrete time Toda lattice and the Askey-Wilson polynomials}, 
Methods and Applications of Analysis {\bf 2}, 369-398 (1995). 

\bi{ode} Odesskii A.V., {\it An analogue of the Sklyanin algebra}, Funkt.
Analiz i Prilozh. {\bf 20}, 78-79 (1986).

\bi{fairlie} Fairlie D.B., {\it Quantum  deformations of $SU(2)$}, 
J. Phys. A: Math. Gen. {\bf 23}, L183-L187 (1990). 

\bi{z} Zhedanov A., {\it Quantum $su_q(2)$ algebra: ``Cartesian" version
and overlaps}, Mod. Phys. Lett. {\bf A7}, 2589-2593 (1992).

\bi{zakh} Zakhariev B.N., {\it Discrete and continuous quantum mechanics.
Exactly solvable models (Lessons of quantum intuition II)},
 Sov. J. Part. Nucl. {\bf{23}}, 603-640 (1992).

\bi{wig} Wiegmann P.B. and Zabrodin A.V., {\it Algebraization of difference
eigenvalue equations related to $U_q(sl_2)$}, Nucl. Phys. {\bf B451},
699-724 (1995).

\bi{zaw} Zhedanov A., {\it Hidden symmetry of the Askey-Wilson polynomials},
Teor. Mat. Fiz. {\bf 89}, 1146-1157 (1991).

\bi{chi} Chihara T.S., {\it An Introduction to Orthogonal Polynomials}
(Gordon and Breach, 1978).

\end {thebibliography}
\end{document}